\documentclass[12pt]{article}

\usepackage[margin=1in]{geometry}  
\usepackage{graphicx}              
\usepackage{amsmath}               
\usepackage{amsfonts}              
\usepackage{amsthm}                
\usepackage{mathtools}
\usepackage{xcolor}
\usepackage{hyperref}

\newtheorem{thm}{Theorem}[section]
\newtheorem{lem}[thm]{Lemma}
\newtheorem{prop}[thm]{Proposition}
\newtheorem{cor}[thm]{Corollary}

\newtheorem{defn}{Definition}

\newtheorem{remark}{Remark}
\newtheorem{fact}{Fact}

\begin{document}

	\nocite{*}
	
	\title{A formula for the $\alpha$-Futaki character}
	
	\author{Kartick Ghosh\thanks{kartick@math.tifr.res.in, School of Mathematics, TIFR Mumbai.}} 
	
	\maketitle
	
	\begin{abstract}
 \'Alvarez-C\'onsul--Garcia-Fernandez--Garc\'ia-Prada introduced the K\"ahler-Yang-Mills equations. They also introduced the $\alpha$-Futaki character, an analog of the Futaki invariant, as an obstruction to the existence of the K\"ahler-Yang-Mills equations. The equations depend on a coupling constant $\alpha$. Solutions of these equations with coupling constant $\alpha>0$ are of utmost importance. In this paper, we provide a formula for the $\alpha$-Futaki character on certain ample line bundles over toric manifolds. We then show that there are no solutions with $\alpha>0$ on certain ample line bundles over certain toric manifolds and compute the value of $\alpha$ if a solution exists. We also relate our result to the existence result of Keller-Friedman in dimension-two.  
	\end{abstract}
	
\section{Introduction}
The K\"ahler-Yang-Mills equations, introduced in (\cite{KYMgeotop}, \cite{mariothesis}), are a system of equations coupling the scalar curvature of K\"ahler metric and a Hermitian-Yang-Mills connections of vector bundles. Suppose $E$ is a holomorphic vector bundle over a K\"ahler manifold $X$. Then for a K\"ahler metric $\omega$ on $X$ and a Hermitian metric $H$ on $E$, the K\"ahler-Yang-Mills equations are 
\begin{equation}
	\label{KYM equations}
	\begin{split}
&\sqrt{-1}{\bigwedge}_{\omega}F_{H}=\lambda Id,\\
&\alpha_{0}S(\omega)-\alpha_{1}{\bigwedge}_{\omega}^{2}tr(F_{H}\wedge F_{H})=c.
	\end{split}
\end{equation}
  Here $\sqrt{-1}{\bigwedge}_{\omega}F_{H}\omega^{n}=n \sqrt{-1}F_{H}\wedge\omega^{n-1}$ and $\lambda,c$ are topological constants. The equations depend on the coupling constants $\alpha_{0}$, $\alpha_{1}$. In this setting, $\frac{\alpha_{1}}{\alpha_{0}}=\alpha$.\\
If $X$ is a compact Riemann surface, the term $tr(F_{H}\wedge F_{H})$ vanishes and hence the system decouples. So for compact Riemann surfaces, the existence problem is equivalent to the  uniformisation theorem and the Narasimhan-Seshadri theorem(\cite{NS}, \cite{DonaldsonNS}).\\
The existence problem in higher dimensions is very difficult because in this case, the system is truly a system of coupled non-linear partial differential equations. Despite this difficulty, there are some existence results: for small $\alpha$ by perturbing a constant scalar curvature K\"ahler metric and a Hermitian-Yang-Mills connection(\cite{KYMgeotop}), applying the dimensional reduction of the K\"ahler-Yang-Mills equations from $\Sigma\times\mathbb{P}^{1}$ to $\Sigma$, where $\Sigma$ is a compact Riemann surface (\cite{KYMcmp},\cite{KYMannalen},\cite{KYMadvances},\cite{KYMBullscimath},\cite{GVreductionbystages}), simultaneous deformation of complex structure of $X$ and $E$(\cite{mario-Tipler}), on line bundles over ruled manifolds using the Calabi ansatz(\cite{keller-friedman}). \\
 \'Alvarez-C\'onsul--Garcia-Fernandez--Garc\'ia-Prada also introduced the $\alpha$-Futaki character $\mathcal{F}$, an analog of the classical Futaki invariant, which is an obstruction to the existence of the K\"ahler-Yang-Mills equations (\ref{KYM equations}). This character does not depend on the choice of a metric $H$ on $E$ and a metric $\omega$ on $X$. The computation of this character is known only in some specific cases(\cite{KYMannalen},\cite{KYMBullscimath}). Also, there is no formula for this character analogous to Mabuchi's formula for the Futaki invariant on toric manifolds(\cite{Mabuchi}) and residue formula for the Futaki invariant(\cite{Tianresidue}, \cite{Futakiresidue}). \\ 
   
 \'Alvarez-C\'onsul--Garcia-Fernandez--Garc\'ia-Prada remarked in \cite{KYMgeotop} that even for line bundles, the existence problem is difficult(page $2807$ in \cite{KYMgeotop}). In this paper, we give a formula of the $\alpha$-Futaki invariant for certain ample line bundles on toric manifolds.\\
 Suppose $(X,\omega)$ is a toric manifold and $L$ is an ample line bundle over it. All the K\"ahler forms will be torus invariant. Take a metric $h$ on $L$ such that its curvature $\beta=\sqrt{-1}F_{h}$ is K\"ahler. We will be concerned with only those ample bundles $L$ such that there exists a metric $\omega'\in [\omega]$ which solves the equation 
 \[n \beta \wedge \omega'^{n-1}=\lambda\omega'^{n}.\]  
This equation is known as the $J$-equation which was introduced by Donaldson(\cite{DonaldsonJ}), Chen(\cite{ChenJ}) and studied on toric manifolds by Yao(\cite{Yi-Yao}) and Collins-Sz\'ekelyhidi(\cite{Gabor-Collins}) . \\
If we calculate the invariant $\mathcal{F}$ with this choice of metrics $h$ and $\omega'$, then we get our main theorem.
\begin{thm}
	\label{main theorem intro}
	Let $P$ and $Q$ be the moment polytopes of $\omega'$ and $\beta$ respectively. Then there is a diffeomorphism $U:P\rightarrow Q$ such that $tr(DU)=\lambda$ and the following holds 
		\[(2\pi)^{-n}\langle \mathcal{F}, \xi_{i}\rangle= \frac{\alpha_{0}}{2} \int_{\partial P}(x_{i}+c_{i})d\sigma+\alpha_{1}\int_{P}(x_{i}+c_{i})(\sum_{m\in M_{2}}m) d\mu.\]
		Here $M_{2}$ is the set of all $(2\times 2)$ principal minor of $DU$, $\xi_{i}$ is a holomorphic vector field on the associated principal $\mathbb{C}^{*}$-bundle of $L$,  $d\mu$ is the Lebesgue measure on $P$, $d\sigma$ is a positive measure on $\partial P$ normalized so that on a face $l_{i}=0$ we have $d\sigma \wedge dl_{i}=\pm d\mu$, $c_{i}$ is a constant such that $\int_{P}(x_{i}+c_{i})d\mu=0$.
\end{thm}
We explain the notations in  section (\ref{section 2}) and section (\ref{section 3}). As an immediate corollary of this theorem, we get the following. 
\begin{cor}
	\label{main cor intro}
		If Solution of K\"ahler-Yang-Mills equations exist with $\alpha_{0}\alpha_{1}\neq 0$ then $\langle \mathcal{F}, \xi_{i}\rangle=0 \ \ \forall i$. Therefore, we get
	\[\frac{\int_{\partial P}(x_{i}+c_{i})d\sigma}{\int_{P}(x_{i}+c_{i})\sum_{m\in M_{2}} m d\mu}=\frac{\int_{\partial P}(x_{j}+c_{j})d\sigma}{\int_{P}(x_{j}+c_{j})\sum_{m\in M_{2}} m d\mu},\]
for all	$ 1\le i,j\le n$.
\end{cor}
 On $\mathbb{P}^2 \# \overline{\mathbb{P}^2}$, we compute the $\alpha$-Futaki character for certain ample line bundles. We prove the following theorem. 
 \begin{thm}
 	\label{dim 2 intro}
 	Suppose $L$ is an ample line bundle on $\mathbb{P}^2 \# \overline{\mathbb{P}^2}$ with $c_{1}(L)=a[H_{2}]-[E_{2}]$ and $[\omega]=b[H_{2}]-[E_{2}]$ satisfying $2\frac{ab-1}{b^{2}-1}>1$. Then there is no solution of the K\"ahler-Yang-Mills equations with $\frac{\alpha_{1}}{\alpha_{0}}>0$. Furthermore, If solution exists then $\frac{\alpha_{1}}{\alpha_{0}}=-\frac{(b^{2}-1)}{(b-a)^{2}}$.
 \end{thm}
We explain the notations in \ref{section 4}. For some of these line bundles, there do exist a solution of the K\"ahler-Yang-Mills equations with the prescribed value of $\frac{\alpha_{1}}{\alpha_{0}}$. The existence result is due to Keller-Friedman(\cite{keller-friedman}).  
\begin{thm}[Keller-Friedman]
	\label{keller-friedman intro}
	There exists metrics $h$ on the line bundle corresponding to the integral cohomology class $(8k+3)[H_{2}]-[E_{2}]$ and $\omega\in 3[H_{2}]-[E_{2}]$ such that they satisfy the K\"ahler-Yang-Mills equations with $\frac{\alpha_{1}}{\alpha_{0}}=-\frac{1}{8k^{2}}$(Here $k\in \mathbb{N}$). 
\end{thm}
We also compute the $\alpha$-Futaki character on some ample line bundles on $\mathbb{P}^3 \# \overline{\mathbb{P}^3}$.
\begin{thm}
	\label{dim 3 intro}
	Suppose $L$ is an ample line bundle on $\mathbb{P}^3 \# \overline{\mathbb{P}^3}$ with $c_{1}(L)=a[H_{3}]-[E_{3}]$ and $[\omega]=b[H_{3}]-[E_{3}]$ satisfying $3\frac{ab^{2}-1}{b^{3}-1}>2$. Then there is no solution of K\"ahler-Yang-Mills equations with $\frac{\alpha_{1}}{\alpha_{0}}>0$. Furthermore, If solution exists then $\frac{\alpha_{1}}{\alpha_{0}}=-\frac{(3b+1)(b-1)^{3}(b^{2}+b+1)}{3b(b-a)^{2}(b^{3}+b^{2}-b+1)}$.
\end{thm}
The author introduced coupled K\"ahler-Einstein and Hermitian-Yang-Mills equations(\cite{kartick}). Suppose $(X, \omega)$ is a compact K\"ahler $n$-manifold with $-\sqrt{-1}\omega$ as the curvature form of the Chern connection of $(K_{X}^{-1},h)$. One can view the metric $h$ as a volume form $\Omega_{h}$.
 Let $\pi: E \rightarrow X$ be a holomorphic vector bundle and  $H$ be a Hermitian metric on the bundle. The coupled K\"ahler-Einstein and Hermitian-Yang-Mills equations are
 \begin{equation}
 	\label{Coupled equations}
 	\begin{split}
 		&\sqrt{-1}{\bigwedge}_{\omega}F_{H}=\lambda Id\\
 		&\frac{\alpha_{0}}{2}\left(\frac{\Omega_{h}}{\int_{X}\Omega_{h}}-\frac{\omega^{n}}{\int_{X}\omega^{n}}\right)-\alpha_{1}\frac{2}{(n-2)!}tr(F_{H}\wedge F_{H})\wedge\omega^{n-2}=\tilde{C}vol_{\omega},
 	\end{split}
 \end{equation}
 where $\alpha_{0},\alpha_{1}$ are coupling constants and $\tilde{C}$, $\lambda$ are topological constants. One can show the following. 
 \begin{lem}[Lemma $3.1.1$ in \cite{kartick}]
Suppose $(X,\omega)$ is a Fano manifold and $f$ is the Ricci potential that is, $Ric(\omega)-\omega=\sqrt{-1}\partial\bar{\partial} f$ normalized so that $\int_{X}e^{f}\omega^{n}=\int_{X}\omega^{n}=V$. Then 
\[\int_{X} \theta_{v}S(\omega)\omega^{n}=V\int_{X}\theta_{v} \left(\frac{\Omega_{h}}{\int_{X}\Omega_{h}}-\frac{\omega^{n}}{\int_{X}\omega^{n}}\right)\]
where $\theta_{v}$ is the Hamiltonian function corresponding to a holomorphic vector field $v$ (that is, $i_{v}\omega=\bar{\partial} \theta_{v}$) and $h$ is a metric on $K_{X}^{-1}$ whose curvature is $\omega$. The function $\theta_{v}$ is normalized by requiring $\int_{X}\theta_{v}\omega^{n}=0$.
  \end{lem}
With the help of this lemma, one can identify the $\alpha$-Futaki character as an obstruction for the existence of coupled K\"ahler-Einstein and Hermitian-Yang-Mills equations albeit when $\omega\in c_{1}(X)$ and $X$ is Fano. The author proved an existence result for line bundles over $\mathbb{P}(\mathcal{O}\oplus\mathcal{O}(-1))\rightarrow \mathbb{P}^{1}$. In this paper, we prove that these line bundles are not ample.\\
Our paper is organized as follows: In section two, we discuss about some basic facts of toric K\"ahler geometry, Futaki invariants, and $J$-equation on toric manifolds. In section three, we recall the $\alpha$-Futaki character and prove theorem \ref{main theorem intro} and corollary \ref{main cor intro}. In section four, we prove theorem \ref{dim 2 intro}, \ref{keller-friedman intro}, \ref{dim 3 intro}.\\
 \textbf{Acknowledgements:} The author would like to thank Prof. Vamsi Pritham Pingali for various useful conversations about the project. The author would also like to thank him for reading a draft of the paper and giving valuable feedbacks. The author would like to thank Prof. Ved V. Datar for a fruitful discussion at IIT Guwahati during the $38$th RMS annual conference. The author would like to thank the organizers of the $38$th RMS annual conference for the invitation. The author was supported by a post-doctoral fellowship from the Tata Institute of Fundamental Research, Mumbai.
 \section{Preliminaries}
 \label{section 2}
 In this section, we briefly discuss some basic notions of toric K\"ahler geometry, the Futaki invariant and  the $J$-equation on toric manifolds. 
 \subsection{K\"ahler geometry of toric manifolds}
\begin{defn}
	A K\"ahler manifold $(X,\omega)$ is called a toric manifold if it admits an effective Hamiltonian action of the torus $(S^{1})^{n}$ which extends to a holomorphic action of the complexified torus $(\mathbb{C}^{*})^{n}$ with a free open dense orbit $\approx (\mathbb{C}^{*})^{n}$.
\end{defn}

The Hamiltonian action induces a moment map $\mu: X \rightarrow \mathbb{R}^{n}$, where the dual of the Lie algebra of the torus is identified with $\mathbb{R}^{n}$. By the Atiyah-Guillemin-Sternberg theorem, the closure of the image is a convex polytope. The size of the polytope is determined by the K\"ahler class while its position is determined by a choice of moment map or normalization of the K\"ahler potential. There is a one-to-one correspondence between toric manifolds and a class of convex polytopes known as Delzant polytope.
\begin{defn}
	A convex polytope is said to be a Delzant polytope if every vertex of the polytope has a neighbourhood which is $SL(n,\mathbb{Z})$ equivalent to $\{y\in \mathbb{R}^{n}| y_{i}\geq 0; i=1,\dots,n\}$.
\end{defn}
A Delzant polytope is said to be integral if all the vertex points are in $\mathbb{Z}^{n}$. A Delzant polytope $P$ is defined by
\begin{equation}
	\label{polytope}
P=\{y\in \mathbb{R}^{n}| l_{i}(y)=\langle y, v_{i}\rangle+\lambda_{i}> 0; i=1,\dots,N\},
\end{equation}
where $v_{i}\in \mathbb{Z}^{n}$ are primitive integral vectors and $\lambda_{i}\in \mathbb{R}$. All the irreducible toric divisors $D_{i}$ are given by the inverse images of the facets $l_{i}=0$ under the moment map. The cohomology class of the K\"ahler form $\omega$ is given by 
\[[\omega]=\sum_{i=1}^{N}\lambda_{i}[D_{i}],\]
where $[D_{i}]$ are Poincar\`e dual of $D_{i}\in H_{2n-2}(X,\mathbb{Z})$. If the Delzant polytope is integral then the K\"ahler class is integral and hence first Chern class of some ample line bundle. \\
Torus invariant K\"ahler metrics take a very special form on toric manifolds. A torus ($(S^{1})^{n}$) invariant metric $\omega$ is of the form $\omega|_{(\mathbb{C}^{*})^{n}}=\sqrt{-1}\partial\bar{\partial}F$; where $F:\mathbb{R}^{n}\rightarrow \mathbb{R}$ is a strictly convex function. The function $F$ is called a K\"ahler potential of $\omega$. Consider the coordinates $x_{i}=\log \lvert z_{i}\rvert^{2}$ and $\theta_{i}$, where $(z_{1}, \dots , z_{n})$ are holomorphic coordinates of $(\mathbb{C}^{*})^{n}$. In this coordinates, the torus invariant metric takes the form
\[\omega|_{(\mathbb{C}^{*})^{n}}=\sqrt{-1}\frac{\partial^{2}F}{\partial x_{i}\partial x_{j}} dx^{i}\wedge d\theta^{j}.\]
Observe that the positivity of the form $\omega$ is equivalent to the strict convexity of the K\"ahler potential. The moment map up to translation is given by $\nabla F$. Let us consider the coordinate $y=\nabla F$. Since $F$ is strictly convex, we can take its Legendre transform. The Legendre transform of $F$ satisfies 
\[G(y)=\langle x,y\rangle +F(x).\]
The function $G$ is called the symplectic potential of $\omega$ and it is a strictly convex function on the moment polytope.
\begin{thm}
	Let $(X,\omega)$ be a toric K\"ahler manifold. Suppose the corresponding moment polytope is defined by (\ref{polytope}). Then every torus invariant metric in the class $[\omega]$ has a symplectic potential of the form
	\[G=\sum_{i=1}^{N}l_{i}(y)\log l_{i}(y)+v,\]
	where $v$ is smooth upto boundary.
\end{thm}
\subsection{The Futaki invariant}
Futaki introduced an invariant as an obstruction to the existence of K\"ahler-Einstein metric(\cite{futaki}). Later Calabi and Futaki generalized the invariants to any compact K\"ahler manifolds. This generalized invariant is an obstruction to the existence of K\"ahler metrics with constant scalar curvature. It is also a character of the Lie algebra of holomorphic vector fields.
\begin{defn}
	The Futaki invariant is given by 
	\[\mathcal{F}ut(v)=\int_{X} \phi_{v} s(\omega)\omega^{n},\]
	where $s(\omega)$ is the scalar curvature of $\omega$, $\phi_{v}$ is a holomorphy potential of $v$(i.e., $ \bar{\partial} \phi_{v}=i_{v}\omega$). Moreover, if $[\omega]=c_{1}(X)$, then
	\[\mathcal{F}ut_{c_{1}}(v)=\int_{X}v(h_{\omega})\omega^{n},\]
	where $h_{\omega}$ satisfies $Ric(\omega)-\omega=\sqrt{-1}\partial \bar{\partial}h_{\omega}$.
\end{defn}
The map $\mathcal{F}ut$ does not depend on the choice of $\omega$.\\

\subsection{ The $J$-equation on toric manifolds}
	The $J$-equation was introduced by Donaldson(\cite{DonaldsonJ}) and Chen(\cite{ChenJ}). Given a K\"ahler form $\beta$ on a compact K\"ahler manifold $X$, it asks for a K\"ahler metric $\omega$ such that 
	\begin{equation}
		\label{J-equation}	
n\beta\wedge \omega^{n-1}=\lambda \omega^{n},
\end{equation}
where $n=dim_{\mathbb{C}}X$ and $\lambda$ is a topological constant. \\
On Toric manifolds, Collins-Sz\'ekelyhidi(\cite{Gabor-Collins}) proved that $J$-equation is solvable iff certain inequalities are met. Yao also studied the $J$-equation on toric manifolds(\cite{Yi-Yao}). \\
Suppose $\beta,\omega$ are two torus invariant K\"ahler metrics on toric manifolds $X$. They induce two moment maps $\mu_{\omega}: X \rightarrow \bar{\mathcal{P}_{\omega}}$ and  $\mu_{\beta}: X \rightarrow \bar{\mathcal{P}_{\beta}}$. Then there exists a transition map 
\begin{equation}
	U : \bar{\mathcal{P}_{\omega}}\rightarrow  \bar{\mathcal{P}_{\beta}}.
\end{equation}
Suppose $u$ is a symplectic potential for $\omega$ and $f$ is a K\"ahler potential for $\beta$, then on the interior of $P$, $U$ is given by
\[U=\nabla f \circ \nabla u.\]
From the work of Yao(page $1586$ in \cite{Yi-Yao}), we have 
\begin{equation}
	\label{Yao formulas}
tr_{\omega}\beta=\frac{n \beta\wedge \omega^{n-1}}{\omega^{n}}=tr(DU|_{\mu_{\omega}}); \ \  \frac{\beta^{n}}{\omega^{n}}=det DU|_{\mu_{\omega}}; \ \  n(n-1)\frac{\beta^{2}\wedge\omega^{n-2}}{\omega^{n}}=\sum_{m\in M_{2}}m,
\end{equation}
where $M_{2}$ is the set of all $(2\times 2)$ principal minors of $DU|_{\mu_{\omega}}$.

	\section{The $\alpha$-Futaki character with toric symmetry}
	\label{section 3}
	In this section, we recall the $\alpha$-Futaki character introduced by Alvarez-Consul--Garcia-Fernandez--Garcia-Prada(see page $2775$ in \cite{KYMgeotop}). \\
	Suppose $\pi: \mathcal{E}\rightarrow X$ is the associated principal $GL(r,\mathbb{C})$-bundle of a holomorphic vector bundle $E$ over $X$, where $r=rank(E)$. Suppose  $I$ is an integrable almost complex structure on $\mathcal{E}$ that preserves the vertical bundle $V\mathcal{E}(=ker(d\pi))$ and $\check{I}$ is the unique integrable almost complex structure on the manifold $X$ determined by $I$. 
	Any vector field $\xi \in Lie Aut(\mathcal{E},I)$ covers a real holomorphic vector field $\check{\xi}$ on $(X,\check{I})$ which can be written as 
	\begin{equation*}
		\check{\xi}=\eta_{\phi_{1}}+\check{I}\eta_{\phi_{2}}+ \beta,
	\end{equation*}
	for any K\"ahler form $\omega$ in some K\"ahler class in $(X,\check{I})$, where $\eta_{\phi_{i}}$ is the Hamiltonian vector field of $\phi_{i}\in C^{\infty}_{0}(X,\omega)$ on $(X,\omega)$ for $i=1,2$, and $\beta$ is the dual of a $1$-form which is harmonic with respect to the K\"ahler metric $\omega(.,\check{I}.)$.
	\begin{equation}
		\label{alpha character}
		\begin{split}
			&\mathcal{F}_{I}:Lie Aut(\mathcal{E}, I)\rightarrow \mathbb{C}\\
			&\langle\mathcal{F}_{I},\xi\rangle = -4\sqrt{-1}\alpha_{1}\int_{X}tr\left(\theta_{H}\xi\wedge(\sqrt{-1}{\bigwedge}_{\omega}F_{H}-\lambda Id)\right)\frac{\omega^{n}}{n!}\\
			&+\int_{X}\phi \left(\alpha_{0}S_{\omega,\check{I}}\omega^{n}+n(n-1)\alpha_{1}tr(\sqrt{-1}F_{H}\wedge \sqrt{-1}F_{H}\wedge\omega^{n-2})\right)\\
			&-\int_{X}\phi (4\alpha_{1}\lambda)tr(\sqrt{-1}F_{H}\wedge\omega^{n-1}),
		\end{split}
	\end{equation}	where $\phi=\phi_{1}+\sqrt{-1}\phi_{2}$. The functional $\mathcal{F}_{I}$ does not depend on the choice of $\omega$ and $H$ and $\mathcal{F}_{I}$ vanishes whenever the coupled equations have a solution. We omit the $I$ in $\mathcal{F}_{I}$ from now onwards. .
\begin{fact}
	\label{fact}
	If $E$ is a line bundle then $\mathcal{E}$ is a principal $\mathbb{C}^{*}$-bundle. It is a fact that any holomorphic vector field on $X$ always lifts to a $\mathbb{C}^{*}$-invariant holomorphic vector field on a principal $\mathbb{C}^{*}$-bundle(see page $2806$ in \cite{KYMgeotop}). We will use this fact in the next subsection. 
\end{fact}
	
\subsection{Formula for certain ample line bundles on toric manifolds}	
Suppose $L$ is an ample line bundle and $\omega$ is a K\"ahler form on a toric manifold $X$. Also assume that there is a metric  $h$ on $L$ such that for its curvature $\sqrt{-1}F_{h}=\beta>0$, there is a metric $\omega'\in [\omega]$ that solves the $J$-equation. We get the formula if we calculate the $\alpha$-Futaki character with $h$ and $\omega'$. We denote by $P,Q$, the moment polytope of $(X,\omega')$ and $(X,\beta)$ respectively. \\
The manifold $X$ has a free open dense subset $(\mathbb{C}^{*})^{n}$. We denote the holomorphic coordinate of $(\mathbb{C}^{*})^{n}$ by $(z_{1},\dots, z_{n})$. The $(\mathbb{C}^{*})^{n}$-invariant holomorphic vector fields $v_{i}=z_{i}\frac{\partial}{\partial z_{i}}$ naturally extend to the toric manifold $X$. Suppose $\xi_{i}$ is a lift of the vector field $v_{i}$ on the associated principal $\mathbb{C}^{*}$ bundle $\mathcal{L}$(using the fact \ref{fact}).  Then we have the following formula.
 \begin{thm}
 	\label{formula theorem}
  There is a diffeomorphism $U:P\rightarrow Q$ such that $tr(DU)=\lambda$ and the following holds 
 \[(2\pi)^{-n}\langle \mathcal{F}, \xi_{i}\rangle= \frac{\alpha_{0}}{2} \int_{\partial P}(x^{i}+c_{i})d\sigma+ \alpha_{1}\int_{P}(x^{i}+c_{i})(\sum_{m\in M_{2}}m) d\mu.\]
 Here $M_{2}$ is the set of all $(2\times 2)$ principal minor of $DU$, $d\mu$ is the Lebesgue measure on $P$, $d\sigma$ is a positive measure on $\partial P$ normalized so that on a face $l_{i}=0$ we have $d\sigma \wedge dl_{i}=\pm d\mu$, and $\int_{P}(x_{i}+c_{i})d\mu=0$.
 \end{thm}
\begin{proof}
	Since for $\sqrt{-1}F_{h}$,  $\omega'$ solves the $J$-equation, we have $\sqrt{-1}\bigwedge_{\omega'}F_{h}=\lambda$. Putting this on $\ref{alpha character}$, we get
	\begin{equation}
		\label{initial form}
	\langle\mathcal{F},\xi_{i}\rangle = 
	\int_{X}\phi_{i} \left(\alpha_{0}S(\omega')\omega'^{n}+\alpha_{1}n(n-1)(\sqrt{-1}F_{h}\wedge \sqrt{-1}F_{h}\wedge\omega'^{n-2})\right),
\end{equation}
	where $i_{v_{i}}\omega'=\bar{\partial}\phi_{i}$ and  $\int_{X}\phi_{i}\omega'^{n}=0$. Now we write down how some terms correspond to similar terms in the polytope $P$. The top form $\omega'^{n}$ corresponds to the measure $(2\pi)^{n}d\mu$ on $P$, where $d\mu$ is the Lebesgue measure on $P$. The function $\phi_{i}\in C^{\infty}_{0}(X,\omega')$ corresponds to the affine function $x^{i}+c_{i}$ with normalisation condition $\int_{P}(x^{i}+c_{i})d\mu=0$. Suppose $u$ is the symplectic potential of $\omega'$, then from Abreu's equation(\cite{Abreu}), we have 
	\[S(\omega')=-\frac{1}{2}\sum_{k,l} \frac{\partial^{2} u^{kl} }{\partial x_{k}\partial x_{l}}.\]
	Hence the first term of \ref{initial form} becomes
	\begin{equation}
			\label{toricfutakiformula}
		\begin{split}
			&\int_{X}\phi_{i}S(\omega')\omega'^{n}\\
			&=\frac{(2\pi)^{n}}{2}\int_{\partial P} (x^{i}+c_{i})d\sigma.
\end{split}
\end{equation}
	We have used Donaldson's formula(\cite{Donaldsonf}, Lemma $4.37$ in \cite{Gaborbook}):
	If $u$ is a symplectic potential and $h: P\rightarrow \mathbb{R}$ is a continuous convex function smooth on the interior, then 
	\[\int_{P} u^{kl}h_{kl}d\mu=\int_{\partial P}hd\sigma-\int_{P}h S(\omega')d\mu,\] 
	where $d\mu$ is the Lebesgue measure on $P$ and $d\sigma$ is a positive measure on the boundary $\partial P$ normalized so that on a face defined by $l_{i}(y)=0$ we have $d\sigma\wedge dl_{i}=\pm d\mu$.

Now putting \ref{Yao formulas} and \ref{toricfutakiformula} in \ref{initial form}, we get
	\[(2\pi)^{-n}\langle \mathcal{F}, \xi_{i}\rangle= \frac{\alpha_{0}}{2} \int_{\partial P}(x^{i}+c_{i})d\sigma+ \alpha_{1}\int_{P}(x^{i}+c_{i})(\sum_{m\in M_{2}}m) d\mu\]

\end{proof}
As an immediate consequence of the theorem, we get
\begin{cor}
	\label{necessary condition}
	If Solution of K\"ahler-Yang-Mills equations exist with $\alpha_{0}\alpha_{1}\neq 0$ then $\langle \mathcal{F}, \xi_{i}\rangle=0 \ \ \forall i$. Therefore, we get
	\[\frac{\int_{\partial P}(x_{i}+c_{i})d\sigma}{\int_{P}(x_{i}+c_{i})\sum_{m\in M_{2}} m d\mu}=\frac{\int_{\partial P}(x_{j}+c_{j})d\sigma}{\int_{P}(x_{j}+c_{j})\sum_{m\in M_{2}} m d\mu},\]
	for all	$ 1\le i,j\le n$.
\end{cor}

\section{Solution of the $J$-equation on $\mathbb{P}^n \# \overline{\mathbb{P}^n}$}
\label{section 4}
Fang and Lai(\cite{Fang-Lai}) studied the $J$-equation on $\mathbb{P}^n \# \overline{\mathbb{P}^n}$ imposing Calabi symmetry on both the metrics. Yao (\cite{Yi-Yao}) cast Fang and Lai's solution in the toric setting as an example of his approach. The Calabi symmetry requires $U(n)$ symmetry instead of $U(1)^{n}$ symmetry as in the toric case. So to make the metrics to have $U(n)$ symmetry, the symplectic potentials need to have more symmetries. \\
We take moment polytopes of $\mathbb{P}^n \# \overline{\mathbb{P}^n}$ in the following form
\[P=\{x\in \mathbb{R}^{n}| x^{i}>0, 1<\sum_{i=1}^{n}x^{i}<b\},\]
where $b$ determines the K\"ahler class. In this case, the K\"ahler class is $[\omega]=b[H_{n}]-[E_{n}]$, where $H_{n}$ is the pull-back of the hyperplane divisor on $\mathbb{P}^{n}$ and $E_{n}$ is the exceptional divisor. Suppose $[\beta]=a[H_{n}]-[E_{n}]$ and the corresponding polytope is $Q$.\\
One takes that the symplectic potentials of $\omega$ and $\beta$ have specific forms
\[u=\sum_{i=1}^{n} x^{i}\log x^{i}-X\log X+h(X), \ \ X=\sum_{i=1}^{n} x^{i}, \]
where $h$ is a strictly convex function defined on $[1,b]$ satisfying 
\[h(X)-(X-1)\log (X-1)-(b-X)\log (b-X)\in C^{\infty}[1,b]\]
and 
\[v=\sum_{i=1}^{n} x^{i}\log x^{i}-X\log X+\theta(X),\]
where $\theta$ is defined on $[1,a]$ satisfying similar condition as $h$. \\
We denote the Legendre transform of $h$ and $\theta$  by $p$ and $\eta$ respectively. We define the map $f(X)=\eta'(h'(X))$.  Then $f: [1,b]\rightarrow [1,a]$ is a diffeomorphism.  The transition map $U:P\rightarrow Q$ in this case is given by 
\[U^{i}(x)=f(X)\frac{x^{i}}{X} \ \ \ \ \ \ \ \ \ \ \ \ \ \ for \ \ 1\le i \le n.\] 
Whenever $n \frac{ab^{n-1}-1}{b^{n}-1}>n-1$, the $J$-equation has a solution given by 
\begin{equation}
	\label{solution of j-equation}
	\tilde{f}(r)=Ar+\frac{B}{r^{n-1}},
\end{equation}
where $A=\frac{ab^{n-1}-1}{b^{n}-1}$ and $B=\frac{b^{n}-ab^{n-1}}{b^{n}-1}$.

Hence we have the following.  
\begin{thm}[\cite{Fang-Lai}, Theorem $2.3$]
	\label{fang-lai theorem}
	On $\mathbb{P}^n \# \overline{\mathbb{P}^n}$, we have two K\"ahler classes $b[H_{n}]-[E_{n}]$, $a[H_{n}]-[E_{n}]$. Suppose $\beta\in a[H_{n}]-[E_{n}]$ is a K\"ahler form. If $n \frac{ab^{n-1}-1}{b^{n}-1}> (n-1)$, then the equation \ref{J-equation} has a solution $\omega\in a[H_{n}]-[E_{n}]$ given by \ref{solution of j-equation} and the corresponding map  $U:P\rightarrow Q$ is given by 
	\[U^{i}(x^{1}, \dots, x^{n})=Ax^{i}+\frac{Bx^{i}}{X^{n}},\]
	where $X=\sum_{i=1}^{n} x^{i}$, $A=\frac{ab^{n-1}-1}{b^{n}-1}$, $B=\frac{b^{n}-ab^{n-1}}{b^{n}-1}$, $P$ moment polytope of $\omega$, $Q$ moment polytope of $\beta$.
\end{thm}

\begin{remark}
	\label{important fact}
	Observe that one can work with K\"ahler classes of the form $a[H_{n}]-b[E_{n}]$ with $a>b>0$ in the theorem \ref{fang-lai theorem} instead of K\"ahler classes of the form $p[H_{n}]-[E_{n}]$ with $p>1$. This will change the boundary values of the corresponding ODE(see Proposition $2.1$ in \cite{Fang-Lai}). 
\end{remark}
\subsection{The $\alpha$-Futaki character on $\mathbb{P}^2 \# \overline{\mathbb{P}^2}$ }
We will work with K\"ahler classes of the following form $p[H_{2}]-[E_{2}]$, where $H_{2}$ is the pull back of the hyperplane divisor from $\mathbb{P}^{2}$, $E_{2}$ is the exceptional divisor, and $p>1$. 

\begin{thm}
	\label{second dimension theorem}
	Suppose $L$ is an ample line bundle on $\mathbb{P}^2 \# \overline{\mathbb{P}^2}$ with $c_{1}(L)=a[H_{2}]-[E_{2}]$ and $[\omega]=b[H_{2}]-[E_{2}]$ satisfying $2\frac{ab-1}{b^{2}-1}>1$. Then there is no solution to the K\"ahler-Yang-Mills equations with $\frac{\alpha_{1}}{\alpha_{0}}>0$. Furthermore, If a solution exists then $\frac{\alpha_{1}}{\alpha_{0}}=-\frac{2b^{2}(b^{2}-1)}{(b^{2}-ab)^{2}}$.
\end{thm}
\begin{proof}	
From \ref{fang-lai theorem}, we know that there is a solution of the $J$-equation and a corresponding map between moment polytopes.
The map $U:P_{2}\rightarrow Q_{2}$ is given by
\[U(x^{1}, x^{2})=\big(Ax^{1}+\frac{Bx^{1}}{(x^{1}+x^{2})^{2}}, Ax^{2}+\frac{Bx^{2}}{(x^{1}+x^{2})^{2}}\big),\]
where $P_{2}=\{(x^{1},x^{2})\in \mathbb{R}^{2}| x^{1}\geq 0, x^{2}\geq 0, 1\leq x^{1}+x^{2}\leq b \}$,  $Q_{2}=\{(x^{1},x^{2})\in \mathbb{R}^{2}| x^{1}\geq 0, x^{2}\geq 0, 1\leq x^{1}+x^{2}\leq a \}$, $A=\frac{ab-1}{b^{2}-1}$, and $B=\frac{b^{2}-ab}{b^{2}-1}$.\\
A simple calculation gives us 
\[ DU=
\begin{bmatrix}
	A+ \frac{B(x^{2}-x^{1})}{(x^{1}+x^{2})^{3}} &  -\frac{2Bx^{1}}{(x^{1}+x^{2})^{3}}\\
	-\frac{2Bx^{2}}{(x^{1}+x^{2})^{3}}  &  	A+ \frac{B(x^{1}-x^{2})}{(x^{1}+x^{2})^{3}}
\end{bmatrix}\]
and 
\[Det(DU)=A^{2}-\frac{B^{2}}{(x^{1}+x^{2})^{4}}. \]
A routine calculation gives 
 \begin{equation}
 	\begin{split}
 		&\int_{P_{2}}d\mu=\frac{b^{2}-1}{2}\\
 		&\int_{\partial P_{2}}d\sigma=3b-1\\
 		&\int_{P_{2}}x^{1}d\mu=\frac{b^{3}-1}{6}=\int_{P_{2}}x^{2}d\mu\\
 		&\int_{\partial P_{2}}x^{1}d\sigma=b^{2}=\int_{\partial P_{2}}x^{2}d\sigma\\
 		&c_{1}=-\frac{\int_{P_{2}}x^{1}d\mu}{\int_{P_{2}}d\mu}=-\frac{b^{2}+b+1}{3(b+1)}=c_{2}\\
 		&\int_{P_{2}}(x^{1}+c_{1})det(DU)dx^{1}dx^{2}=B^{2}\frac{(b-1)^{3}}{6b^{2}}=\frac{(b^{2}-ab)^{2}(b-1)}{6b^{2}(b+1)^{2}}
 	\end{split}
 \end{equation}
Now if there is a solution to the K\"ahler-Yang-Mills equations then the $\alpha$-Futaki character vanishes. Hence using theorem \ref{formula theorem}, we get $\frac{\alpha_{1}}{\alpha_{0}}=-\frac{(b^{2}-1)}{(b-a)^{2}}$. \\
Hence there is no solution to the K\"ahler-Yang-Mills equations for $\frac{\alpha_{1}}{\alpha_{0}}>0$ and furthermore if there is a solution then $\frac{\alpha_{1}}{\alpha_{0}}=-\frac{(b^{2}-1)}{(b-a)^{2}}$.
\end{proof}
\begin{remark}
	If we use K\"ahler classes of the form $a[H]-b[E]$ with $a>b>0$ as stated in \ref{important fact},  we will get similar results, that is, $\frac{\alpha_{1}}{\alpha_{0}}<0$.
\end{remark}
 Keller and Friedman constructed solution of K\"ahler-Yang-Mills equations on some line bundles over $\mathbb{P}^2 \# \overline{\mathbb{P}^2}$.
 \begin{thm}[Keller-Friedman]
 	There exists metrics $h$ on the line bundle corresponding to the integral cohomology class $(8k+3)[H_{2}]-[E_{2}]$ and $\omega\in 3[H_{2}]-[E_{2}]$ such that they satisfy the K\"ahler-Yang-Mills equations with $\frac{\alpha_{1}}{\alpha_{0}}=-\frac{1}{8k^{2}}$(Here $k\in \mathbb{N}$). 
 \end{thm}
\begin{proof}
In \cite{keller-friedman}, Keller and Friedman proved that:\\ Suppose $M=\mathbb{P}(\mathcal{O}\oplus \mathcal{L})\rightarrow \Sigma$ is a ruled surface with $\mathcal{L}$ of degree $k\in \mathbb{Z}^{*}_{+}$. Fix $k'\in \mathbb{Z}^{*}_{+}, k_{1}\in \mathbb{Z}, k_{2}\in \mathbb{Z}^{*}$. Consider the integral classes $L=2[E_{0}]+k'[C]$ and $E=2(k_{1}-k_{2})[E_{0}]+(k_{1}k'+k_{2}(2k+k'))[C]$. Then there exists a solution of K\"ahler-Yang-Mills equations with $\frac{\alpha_{1}}{\alpha_{0}}=-\frac{(2-s_{\Sigma})k+2k'}{8k_{2}^{2}(k+k')}$, where $E_{0}$ is the zero section, $C$ denotes fiber, $s_{\Sigma}=\frac{2(1-h)}{k}$, $h=$genus$(\Sigma)$.\\

	Since $\mathbb{P}^2 \# \overline{\mathbb{P}^2}=\mathbb{P}(\mathcal{O}\oplus \mathcal{O}(1))\rightarrow \mathbb{P}^{1}$, we take $\Sigma=\mathbb{P}^{1}$ and $k=1$. The relations between $[E_{2}],[H_{2}],$ and $C, E_{0}$ are $[C]=[H_{2}]-[E_{2}], [E_{0}]=[H_{2}]$. Hence for $L$ to be of the form $a[H_{2}]-[E_{2}]$, we must have $k'=1$. Now $E$ becomes $2(k_{1}-k_{2})[E_{0}]+(k_{1}+3k_{2})[C]=(3k_{1}+k_{2})[H_{2}]-(k_{1}+3k_{2})[E_{2}]$. Now putting $k_{2}=-k$ and $k_{1}=1+3k$ with $k\in \mathbb{N}$, we get $E=(8k+3)[H_{2}]-[E_{2}]$. Finally, we compute $\frac{\alpha_{1}}{\alpha_{0}}=-\frac{1}{8k^{2}}$. Observe that our theorem \ref{second dimension theorem} also asserts that $\frac{\alpha_{1}}{\alpha_{0}}$ has to be $-\frac{3^{2}-1}{(3-8k-3)^{2}}=-\frac{1}{8k^{2}}$.
\end{proof}
\begin{remark}
	The above analysis also indicates that why Keller-Friedman could not find solution of K\"ahler-Yang-Mills equations with $\frac{\alpha_{1}}{\alpha_{0}}>0$.
\end{remark}

\subsection{The $\alpha$-Futaki character on $\mathbb{P}^3 \# \overline{\mathbb{P}^3}$ }

We will work with K\"ahler classes of the form $p[H_{3}]-[E_{3}]$, where $H_{3}$ is the pullback of the hyperplane divisor from $\mathbb{P}^{3}$ and $E_{3}$ is the exceptional divisor, $p>1$. 

\begin{thm}
	Suppose $L$ is an ample line bundle on $\mathbb{P}^3 \# \overline{\mathbb{P}^3}$ with $c_{1}(L)=a[H_{3}]-[E_{3}]$ and $[\omega]=b[H_{3}]-[E_{3}]$ satisfying $3\frac{ab^{2}-1}{b^{3}-1}>2$. Then there is no solution to K\"ahler-Yang-Mills equations with $\frac{\alpha_{1}}{\alpha_{0}}>0$. Furthermore, If a solution exists then $\frac{\alpha_{1}}{\alpha_{0}}=-\frac{(3b+1)(b-1)^{3}(b^{2}+b+1)}{3b(b-a)^{2}(b^{3}+b^{2}-b+1)}$.
\end{thm}
\begin{proof}
	From \ref{fang-lai theorem}, we know that there is a solution of the $J$-equation and a corresponding map between moment polytopes. So we have a map $U:P_{3}\rightarrow Q_{3}$ given by
\[U(x^{1}, x^{2}, x^{3})=\big(Ax^{1}+\frac{Bx^{1}}{(x^{1}+x^{2}+x^{3})^{3}}, Ax^{2}+\frac{Bx^{2}}{(x^{1}+x^{2}+x^{3})^{3}}, Ax^{3}+\frac{Bx^{3}}{(x^{1}+x^{2}+x^{3})^{3}}\big),\]
where $P_{3}=\{(x^{1},x^{2},x^{3})\in \mathbb{R}^{3}| x^{i}\geq 0 \ i=1,2,3; 1\leq x^{1}+x^{2}+x^{3}\leq b \}$,  $Q_{3}=\{(x^{1},x^{2},x^{3})\in \mathbb{R}^{3}| x^{i}\geq 0\  i=1,2,3; 1\leq x^{1}+x^{2}+x^{3}\leq a \}$, $A=\frac{ab^{2}-1}{b^{3}-1}$, and $B=\frac{b^{3}-ab^{2}}{b^{3}-1}$.\\
A simple calculation yields
\[DU=
\begin{bmatrix}
	A-\frac{3Bx^{1}}{X^{4}}+\frac{B}{X^{3}} & -\frac{3Bx^{1}}{X^{4}} & -\frac{3Bx^{1}}{X^{4}}\\
-\frac{3Bx^{2}}{X^{4}} &	A-\frac{3Bx^{2}}{X^{4}}+\frac{B}{X^{3}} &  -\frac{3Bx^{2}}{X^{4}}\\
 -\frac{3Bx^{3}}{X^{4}} & -\frac{3Bx^{3}}{X^{4}} &	A-\frac{3Bx^{3}}{X^{4}}+\frac{B}{X^{3}} \\
\end{bmatrix},\]
where $X=x^{1}+x^{2}+x^{3}$. Suppose $M_{1}, M_{2}, M_{3}$ are all the $2\times 2 $ principal minor of $DU$. Then 
\[M_{1}+M_{2}+M_{3}=3(A^{2}-\frac{B^{2}}{X^{6}}).\]
A routine calculation yields 
\begin{equation}
	\begin{split}
		& \int_{P_{3}}d\mu=\frac{b^{3}-1}{6}\\
		& \int_{\partial P_{3}} d\sigma=2b^{2}-1\\
		& \int_{P_{3}} x^{i}d\mu=\frac{b^{4}-1}{24},\ \ i=1,2,3\\
		& \int_{\partial P_{3}} x^{i}d\sigma=\frac{3b^{3}-1}{6}, \ \  i=1,2,3\\
		& \int_{P_{3}} (x^{i}+c_{i})3(A^{2}-\frac{B^{2}}{X^{6}})d\mu=\frac{B^{2}(b^{4}-2b^{3}+2b-1)}{8b^{3}},\ \ i=1,2,3
	\end{split}
\end{equation}
where $c_{i}=-\frac{(b^{2}+1)(b+1)}{4(b^{2}+b+1)}$.\\
Now if there is a solution to the K\"ahler-Yang-Mills equations then the $\alpha$-Futaki character vanishes. Hence using theorem \ref{formula theorem}, we get $\frac{\alpha_{1}}{\alpha_{0}}=-\frac{(3b+1)(b-1)^{3}(b^{2}+b+1)}{3b(b-a)^{2}(b^{3}+b^{2}-b+1)}$.
Hence there is no solution of the K\"ahler-Yang-Mills equations for $\frac{\alpha_{1}}{\alpha_{0}}>0$ and furthermore, if there is a solution then 
$\frac{\alpha_{1}}{\alpha_{0}}=-\frac{(3b+1)(b-1)^{3}(b^{2}+b+1)}{3b(b-a)^{2}(b^{3}+b^{2}-b+1)}$.
\end{proof}

Finally, we prove that the line bundles that we constructed in \cite{kartick} over\\ $\mathbb{P}(\mathcal{O}\oplus\mathcal{O}(-1))\rightarrow \mathbb{P}^{1}$ are not ample.
\begin{prop}
	Suppose $X=\mathbb{P}(\mathcal{O}\oplus \mathcal{O}(-1))\rightarrow \mathbb{P}^{1}$ and $L$ is the line bundle corresponding to the second integral cohomology class $\delta=\frac{a\omega+b\gamma}{2\pi}$, where $\omega$ is a K\"ahler metric, $\gamma$ is a tracefree(with respect to $\omega$) $(1,1)$ form,  $a=\frac{m_{1}+m_{2}\log 3}{2+3\log 3}$, and  $b=\frac{2m_{2}-3m_{1}}{2+3\log 3}$ for $(m_{1},m_{2})\ne (0,0)$. These line bundles are not ample.
\end{prop}
\begin{proof}
	First we recall the expression for $\omega$ and $\gamma$(for details see \cite{kartick}). The expressions are as follows
	\[\omega=(1+\tau)\omega_{FS}+\phi(\tau)\frac{\sqrt{-1}dw\wedge d\bar{w}}{\lvert w\rvert^{2}},\]
	 \[\gamma=\omega_{FS}-\frac{\phi(\tau)}{1+\tau}\frac{\sqrt{-1} dw\wedge d\bar{w}}{\lvert w\rvert^{2}},\]
	 where $\phi(\tau):[0,2]\rightarrow \mathbb{R}$ is a function, $\tau:X\rightarrow [0,2]$ is a map, $w$ fiber coordinate, and $\omega_{FS}$ is the Fubini-Study metric on $X$. \\
	 Let's assume that the line bundles are ample. Then from Nakai-Moishezon criterion, we get 
	 \[\int_{\mathbb{P}^{1}} \delta>0, \ \  \int_{C}\delta>0, \ \  \int_{X} \delta^{2}>0, \]
	 where $C$ denotes fiber.\\
	 A routine calculations yields
	 \[a+b>0, \ \  2a-b\log 3>0, \ \ b^{2}\log 3 > 4a^{2}.  \]
	 Writing these inequalities in terms of $m_{1},m_{2}$, we get 
	 \[m_{1}<(1+\log \sqrt{3})m_{2}, \ \ m_{1}>0, \ \ (9\log 3-4)m_{1}^{2}>20m_{1}m_{2}\log 3.\]
	 Now 
	 \begin{equation}
	 	\begin{split}
	 		&m_{1}<(1+\log \sqrt{3})m_{2}<\frac{9\log 3-4}{20 \log 3}m_{1}\\
	 		&20\log 3< 9\log 3-4\\
	 		& 4+11\log 3<0. 
	 	\end{split}
	 \end{equation}
	 A contradiction.
\end{proof}
		\bibliographystyle{plain}
			
			\bibliography{paper}

		\end{document}